\documentclass[12pt]{amsart}

\usepackage[paperheight=11in,paperwidth=8.5in,left=1in,top=1.25in,right=1in,bottom=1.25in,]{geometry}
\usepackage[linktocpage=false,colorlinks=true,linkcolor=black,citecolor=black,
filecolor=black,urlcolor=black,pagebackref=false,pdfstartpage={1},
pdfstartview={FitH},pdftitle={},pdfauthor={},pdfsubject={},pdfcreator={},
pdfproducer={},pdfkeywords={}]{hyperref}
\usepackage{afterpage}
\usepackage{amsfonts}
\usepackage{amsmath}
\allowdisplaybreaks[4]
\usepackage{amssymb}
\usepackage{amsthm}
\usepackage{array}
\usepackage{bbm}
\usepackage[justification=centering]{caption}
\usepackage[capitalize,nameinlink,noabbrev,nosort]{cleveref}
\usepackage{enumerate}
\usepackage{enumitem}
\usepackage{float}
\restylefloat{figure}
\usepackage{graphicx}
\usepackage{mathrsfs}
\usepackage[framemethod=tikz,ntheorem,xcolor]{mdframed}
\usepackage{parcolumns}
\usepackage{tabularx}
\usepackage{tikz}\pdfpageattr{/Group <</S /Transparency /I true /CS /DeviceRGB>>}
\usetikzlibrary{arrows,calc,intersections,matrix,positioning,through}
\usepackage{tikz-cd}
\tikzset{commutative diagrams/.cd,every label/.append style = {font = \normalsize}}
\usepackage[all]{xy}
\usepackage[aligntableaux=center]{ytableau}

\numberwithin{equation}{section}

\newtheorem{thm}[equation]{Theorem}
\AfterEndEnvironment{thm}{\noindent\ignorespaces}

\AfterEndEnvironment{conj}{\noindent\ignorespaces}
\newtheorem{cor}[equation]{Corollary}
\AfterEndEnvironment{cor}{\noindent\ignorespaces}
\newtheorem{lem}[equation]{Lemma}
\AfterEndEnvironment{lem}{\noindent\ignorespaces}
\newtheorem{prob}[equation]{Problem}
\AfterEndEnvironment{prob}{\noindent\ignorespaces}
\newtheorem{prop}[equation]{Proposition}
\AfterEndEnvironment{prop}{\noindent\ignorespaces}

\theoremstyle{definition}
\newtheorem{defn}[equation]{Definition}
\AfterEndEnvironment{defn}{\noindent\ignorespaces}
\newtheorem*{pf_no_qed}{Proof}
\newenvironment{pf}[1][]{\begin{pf_no_qed}[#1]\pushQED{\qed}}{\popQED\end{pf_no_qed}}
\AfterEndEnvironment{pf}{\noindent\ignorespaces}
\newtheorem{eg_no_qed}[equation]{Example}
\newenvironment{eg}[1][]{\begin{eg_no_qed}[#1]\pushQED{\qed}}{\popQED\end{eg_no_qed}}
\AfterEndEnvironment{eg}{\noindent\ignorespaces}
\newtheorem{rmk}[equation]{Remark}
\AfterEndEnvironment{rmk}{\noindent\ignorespaces}

\theoremstyle{remark}

\AfterEndEnvironment{claim}{\noindent\ignorespaces}
\newtheorem*{claimpf_no_qed}{Proof of Claim}

\AfterEndEnvironment{claimpf}{\noindent\ignorespaces}

\newcommand{\cd}{\cdots\!}

\DeclareMathOperator{\GL}{GL}
\DeclareMathOperator{\gl}{\mathfrak{gl}}
\DeclareMathOperator{\Gr}{Gr}

\DeclareMathOperator{\pr}{pr}

\newcommand{\rf}[1]{\hyperref[#1]{(\ref*{#1})}}

\DeclareMathOperator{\spn}{span}

\DeclareMathOperator{\tr}{tr}

\title{Moment curves and cyclic symmetry for positive Grassmannians}

\author{Steven N. Karp}
\address{Department of Mathematics, University of Michigan}
\email{\href{mailto:snkarp@umich.edu}{snkarp@umich.edu}}

\thanks{This work was supported in part by NSF grant DMS-1600447.}

\begin{document}

\begin{abstract}
We show that for each $k$ and $n$, the cyclic shift map on the Grassmannian $\Gr_{k,n}(\mathbb{C})$ has exactly $\binom{n}{k}$ fixed points. There is a unique totally nonnegative fixed point, given by taking $n$ equally spaced points on the trigonometric moment curve (if $k$ is odd) or the symmetric moment curve (if $k$ is even). We introduce a parameter $q\in\mathbb{C}^\times$, and show that the fixed points of a $q$-deformation of the cyclic shift map are precisely the critical points of the mirror-symmetric superpotential $\mathcal{F}_q$ on $\Gr_{k,n}(\mathbb{C})$. This follows from results of Rietsch about the quantum cohomology ring of $\Gr_{k,n}(\mathbb{C})$. We survey many other diverse contexts which feature moment curves and the cyclic shift map.
\end{abstract}

\maketitle

\section{Introduction}\label{sec_introduction}

\noindent The {\itshape (complex) Grassmannian $\Gr_{k,n}(\mathbb{C})$} is the set of $k$-dimensional subspaces of $\mathbb{C}^n$. We represent $V\in\Gr_{k,n}(\mathbb{C})$ by a $k\times n$ matrix $A$, unique up to row operations, whose rows form a basis of $V$. For $k$-subsets $I\subseteq\{1, \cd, n\}$ we let $\Delta_I(V)$ be the $k\times k$ minor of $A$ located in columns $I$. The minors $\Delta_I(V)$ do not depend on our choice of $A$ up to a common nonzero scalar, and give projective coordinates on $\Gr_{k,n}(\mathbb{C})$, called {\itshape Pl\"{u}cker coordinates}. We call $V$ {\itshape totally nonnegative} if $\Delta_I(V)\ge 0$ for all $I$ (for some choice of the common scalar), and {\itshape totally positive} if $\Delta_I(V) > 0$ for all $I$. For example, the span of $(1,0,0,-1)$ and $(-1,2,1,3)$ is a totally nonnegative element of $\Gr_{2,4}(\mathbb{C})$, but it is not totally positive. The set of all totally nonnegative $V\in\Gr_{k,n}(\mathbb{C})$ forms the {\itshape totally nonnegative Grassmannian} $\Gr_{k,n}^{\ge 0}$.

The totally nonnegative Grassmannian has a stratification into cells \cite{rietsch99} exhibiting remarkable combinatorial properties \cite{postnikov}, one of which is {\itshape cyclic symmetry}. For fixed $k$ and $n$, we define the {\itshape (left) cyclic shift map $\sigma\in\GL_n(\mathbb{C})$} by
$$
\sigma(v) := (v_2, v_3, \cd, v_n, (-1)^{k-1}v_1) \quad \text{ for } v = (v_1, \cd, v_n)\in\mathbb{C}^n.
$$
Given $V\in\Gr_{k,n}(\mathbb{C})$, we denote by $\sigma(V)$ the subspace $\{\sigma(v) : v\in V\}\in\Gr_{k,n}(\mathbb{C})$, so $\sigma$ acts on Pl\"{u}cker coordinates by rotating the index set $\{1, \cd, n\}$. Hence $\sigma$ is an automorphism of $\Gr_{k,n}(\mathbb{C})$ of order $n$, which restricts to an automorphism of $\Gr_{k,n}^{\ge 0}$. We show the following:
\begin{thm}\label{unique_fixed_point}
The cyclic shift map $\sigma$ on $\Gr_{k,n}(\mathbb{C})$ has exactly $\binom{n}{k}$ fixed points, namely, $\spn\{(1, z_j, z_j^2, \cd, z_j^{n-1}) : 1 \le j \le k\}$ for any $k$ distinct $n$th roots $z_1, \cd, z_k\in\mathbb{C}$ of $(-1)^{k-1}$. There is a unique totally nonnegative fixed point $V_0$, obtained by taking $z_1, \cd, z_k$ to be the $k$ roots closest to $1$ on the unit circle.
\end{thm}

We give an elementary proof of this result, deducing the uniqueness of the totally nonnegative fixed point from the following characterization of $\Gr_{k,n}^{\ge 0}$ in terms of sign variation:
\begin{thm}[Gantmakher and Krein {\cite[Theorem V.3]{gantmakher_krein_50}}]\label{gantmakher_krein}
Suppose that all Pl\"{u}cker coordinates of $V\in\Gr_{k,n}(\mathbb{C})$ are real (up to a common scalar). Then $V$ is totally nonnegative if and only if each real vector in the subspace $V$ changes sign at most $k-1$ times, when viewed as a sequence of $n$ numbers and ignoring any zeros.
\end{thm}
This is a Grassmannian analogue of the well-known result of Schoenberg \cite{schoenberg30}, that an injective linear map $A:\mathbb{R}^k\to\mathbb{R}^n$ diminishes sign variation if and only if for each $j$, the nonzero $j\times j$ minors of $A$ have the same sign.

There is an elegant way to describe the unique totally nonnegative fixed point $V_0\in\Gr_{k,n}^{\ge 0}$ in terms of certain real curves. Define $f_k:\mathbb{R}\to\mathbb{R}^k$ by
\begin{align}\label{f_defn}
f_k(\theta) := \begin{cases}
\scalebox{0.96}{$\left(1, \cos\!\left(\theta\right)\!, \sin\!\left(\theta\right)\!, \cos\!\left(2\theta\right)\!, \sin\!\left(2\theta\right)\!, \cd, \cos\!\left(\frac{k-1}{2}\theta\right)\!, \sin\!\left(\frac{k-1}{2}\theta\right)\right)$}, & \text{if $k$ is odd,}\\[2pt]
\scalebox{0.96}{$\left(\cos\!\left(\frac{1}{2}\theta\right)\!, \sin\!\left(\frac{1}{2}\theta\right)\!, \cos\!\left(\frac{3}{2}\theta\right)\!, \sin\!\left(\frac{3}{2}\theta\right)\!, \cd, \cos\!\left(\frac{k-1}{2}\theta\right)\!, \sin\!\left(\frac{k-1}{2}\theta\right)\right)$}, & \text{if $k$ is even.}\end{cases}
\end{align}
Note that $f_k(\theta+2\pi) = (-1)^{k-1}f_k(\theta)$. For odd $k$, the curve in $\mathbb{R}^{k-1}$ formed from $f_k$ by deleting the first component is the {\itshape trigonometric moment curve}, and for even $k$, the curve $f_k$ is the {\itshape symmetric moment curve}. These curves have a rich history, which we discuss in \cref{sec_connections}. The fixed point $V_0$ is represented by any $k\times n$ matrix whose columns are $f_k(\theta_1), \cd, f_k(\theta_n)$, such that the points $\theta_1 < \theta_2 < \cdots < \theta_n < \theta_1 + 2\pi$ are equally spaced on the real line, i.e.\ $\theta_{j+1} - \theta_j = \frac{2\pi}{n}$ for $1 \le j \le n-1$. We observe that by \cref{unique_fixed_point}, for fixed $k$ we can recover the curve $f_k$ (up to an automorphism of $\mathbb{R}^k$) from the cyclic fixed points $V_0\in\Gr_{k,n}^{\ge 0}$, by taking a sort of limit as $n\to\infty$.

We also have the following explicit formula for the Pl\"{u}cker coordinates of $V_0\in\Gr_{k,n}^{\ge 0}$:
\begin{align}\label{V_Plueckers}
\Delta_I(V_0) = \prod_{1 \le r < s \le k}\sin\!\left(\textstyle\frac{i_s - i_r}{n}\pi\right) \quad \text{ for all $k$-subsets $I = \{i_1 < \cdots < i_k\}\subseteq\{1, \cd, n\}$}.
\end{align}
Note that this immediately implies $V_0$ is totally positive, since $\sin(\theta) > 0$ for $0 < \theta < \pi$. We prove \cref{unique_fixed_point} and the stated properties of $V_0$ in \cref{sec_cyclic_shift}.
\begin{eg}\label{intro_eg}
Let $k := 2$, $n := 4$. Then $f_k$ parametrizes the unit circle in the plane, and the fixed point $V_0$ is represented by the matrix $A := \scalebox{0.9}{$\begin{bmatrix}
1 & \frac{1}{\sqrt{2}} & 0 & -\frac{1}{\sqrt{2}} \\[3pt]
0 & \frac{1}{\sqrt{2}} & 1 & \frac{1}{\sqrt{2}}
\end{bmatrix}$}$, whose columns correspond to four consecutive points on the regular unit octagon:
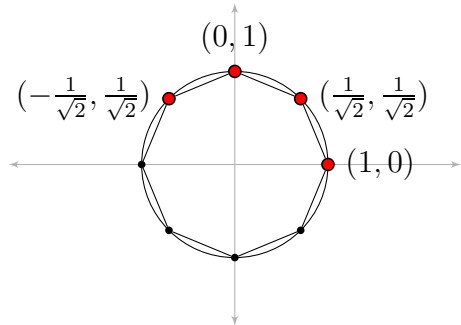
\begin{align*}
\qquad\begin{tikzpicture}[baseline=(current bounding box.center),scale=0.86]
\tikzstyle{vertex}=[inner sep=0,minimum size=0.8mm,circle,draw=black,fill=black,semithick]
\tikzstyle{selected}=[inner sep=0,minimum size=1.6mm,circle,draw=black,fill=red,semithick]
\pgfmathsetmacro{\radius}{1.44};
\path[latex'-latex',color=black!30](-3.5,0)edge(3.5,0) (0,-2.5)edge(0,2.5);
\draw(0,0)circle[radius=\radius];
\node[selected][label=right:{$(1,0)$}](b1)at(0:\radius){};
\node[selected][label=right:{$(\frac{1}{\sqrt{2}},\frac{1}{\sqrt{2}})$}](b2)at(45:\radius){};
\node[selected][label=above:{$(0,1)$}](b3)at(90:\radius){};
\node[selected][label=left:{$(-\frac{1}{\sqrt{2}},\frac{1}{\sqrt{2}})$}](b4)at(135:\radius){};
\foreach \x in {5,...,8}{
\node[vertex](b\x)at(\x*45-45:\radius){};}
\path (b1)edge(b2) (b2)edge(b3) (b3)edge(b4) (b4)edge(b5) (b5)edge(b6) (b6)edge(b7) (b7)edge(b8) (b8)edge(b1);
\end{tikzpicture}\qquad.
\end{align*}
We see that $V_0$ is indeed fixed by $\sigma$, since $\sigma(A)$ is $\scalebox{0.9}{$\begin{bmatrix}
\frac{1}{\sqrt{2}} & 0 & -\frac{1}{\sqrt{2}} & -1 \\[3pt]
\frac{1}{\sqrt{2}} & 1 & \frac{1}{\sqrt{2}} & 0
\end{bmatrix}$}$, which (although not equal to $A$ as a $2\times 4$ matrix) has the same $2\times 2$ minors as $A$.

Alternatively, we can represent $V_0$ by the matrix $A' := \begin{bmatrix}1 & \zeta & \zeta^2 & \zeta^3 \\ 1 & \zeta^{-1} & \zeta^{-2} & \zeta^{-3}\end{bmatrix}$, where $\zeta := e^{\frac{i\pi}{4}}$ and $\zeta^{-1}$ are the two fourth roots of $-1$ closest to $1$ on the unit circle. We can check that the $2\times 2$ minors of $A'$ and $A$ are equal, up to the common scalar $i/2$.
\end{eg}

In \cref{sec_quantum}, we discuss a remarkable connection of the cyclic shift map to the {\itshape quantum cohomology ring} of $\Gr_{k,n}(\mathbb{C})$, which is a new manifestation of the `cyclic hidden symmetry' (see \cite[Section 6.2]{postnikov_05}). Namely, Rietsch \cite{rietsch01} has proved that this ring is isomorphic to the coordinate ring of a certain variety $\mathcal{Y}_{k,n}$, which had been discovered by Peterson in unpublished work. There is a natural embedding of $\mathcal{Y}_{k,n}$ into $\Gr_{k,n}(\mathbb{C})$, and we show that its image is precisely the union of the fixed points of the maps $\sigma_q$ over $q\in\mathbb{C}$, where $\sigma_q$ is a $q$-deformation of the cyclic shift map $\sigma$. It then follows from work of Rietsch that the fixed points of $\sigma_q$ are precisely the critical points of the mirror-symmetric {\itshape superpotential} $\mathcal{F}_q$ of $\Gr_{k,n}(\mathbb{C})$.\footnote{Linhui Shen and Peng Zhou have independently established the same result.} Intriguingly, describing the totally nonnegative part of $\mathcal{Y}_{k,n}$ is equivalent to establishing certain inequalities among Schur polynomials evaluated at roots of unity. Rietsch proved these inequalities using orthogonality relations; we show how to derive them from \cref{gantmakher_krein}.

In \cref{sec_connections}, we survey several topics related to the curves $f_k$ or the cyclic shift map $\sigma$, including cyclic and bicyclic polytopes, isoperimetric inequalities for convex curves, separable quantum states, arrangements of equal minors, the twist map, promotion on rectangular tableaux, birational rowmotion, control theory, and the topology of $\Gr_{k,n}^{\ge 0}$.

~\\\noindent{\bfseries Acknowledgements.} I am grateful to Lauren Williams for her encouragement and many helpful suggestions. I also thank Nima Arkani-Hamed, Gabriel Frieden, Pavel Galashin, Alexander Garver, Michael Gekhtman, Rachel Karpman, Thomas Lam, Alexander Postnikov, Vic Reiner, Brendon Rhoades, Konstanze Rietsch, Linhui Shen, and Franco Vargas Pallete for enlightening discussions, and the reviewers for their valuable feedback.

\section{Fixed points of the cyclic shift map}\label{sec_cyclic_shift}

\noindent In this section we identify all fixed points of the cyclic shift map $\sigma$ on $\Gr_{k,n}(\mathbb{C})$.
\begin{defn}\label{defn_V(S)}
Given a $k$-subset $S\subseteq\mathbb{C}$, define $V_S\in\Gr_{k,n}(\mathbb{C})$ as the subspace with basis $\{(1, z, \cd, z^{n-1}) : z\in S\}$, which has dimension $k$ by Vandermonde's determinantal identity. We also denote the $k$-subset $\left\{e^{-i\frac{(k-1)\pi}{n}}, e^{-i\frac{(k-3)\pi}{n}}, \cd, e^{i\frac{(k-3)\pi}{n}}, e^{i\frac{(k-1)\pi}{n}}\right\}$ of $n$th roots of $(-1)^{k-1}$ closest to $1$ on the unit circle by $S_0$. We let $V_0 := V_{S_0}$.
\end{defn}

First we show that we can uniquely recover $S$ from $V_S$.
\begin{lem}\label{recover_S}
Let $S$ and $S'$ be $k$-subsets of $\mathbb{C}$ such that $V_S = V_{S'}$. Then $S = S'$.
\end{lem}

\begin{pf}
It suffices to show that $S'\subseteq S$. Suppose otherwise that there exists $z'\in S'\setminus S$. Then $\{(1, z, \cd, z^{n-1}) : z\in S\cup\{z'\}\}$ is contained in $V_S$, and it is linearly independent by Vandermonde's determinantal identity. This implies $\dim(V_S) \ge k+1$, a contradiction.
\end{pf}

Now we establish the properties of $V_0$ stated in \cref{sec_introduction}.
\begin{lem}\label{vandermonde}\cite{scott}
For $k\ge 0$ and $\theta_1, \cd, \theta_k\in\mathbb{R}$, we have
$$
\det(f_k(\theta_1), \cd, f_k(\theta_k)) = 2^{\lfloor (k-1)^2/2\rfloor}\prod_{1 \le r < s \le k}\sin\!\left(\textstyle\frac{\theta_s - \theta_r}{2}\right).
$$
\end{lem}
Scott \cite{scott} stated this result for odd $k$, and outlined a proof, which also handles the case of even $k$. We give the complete proof here.
\begin{pf}
Define $g_k:\mathbb{R}\to\mathbb{C}^k$ by $g_k(\theta) := (e^{-\frac{k-1}{2}i\theta}, e^{-\frac{k-3}{2}i\theta}, \cd, e^{\frac{k-3}{2}i\theta}, e^{\frac{k-1}{2}i\theta})$. We regard $f_k(\theta)$ and $g_k(\theta)$ as column vectors in $\mathbb{C}^k$. Using
\begin{align}\label{exp_to_trig}
\begin{bmatrix}
1 & -i \\
1 & i
\end{bmatrix}
\begin{bmatrix}
\cos(\theta) \\
\sin(\theta)
\end{bmatrix}
=
\begin{bmatrix}
e^{-i\theta} \\
e^{i\theta}
\end{bmatrix},
\end{align}
we can get from the matrix $[f_k(\theta_1) \,|\, \cdots \,|\, f_k(\theta_k)]$ to $[g_k(\theta_1) \,|\, \cdots \,|\, g_k(\theta_k)]$, introducing a factor of $\pm(2i)^{-\lfloor k/2\rfloor}$ when we take determinants. (The sign $\pm$ is $-$ if $k$ and $\frac{k-1}{2}$ are both odd, and $+$ otherwise.) When we scale each column of $[g_k(\theta_1) \,|\, \cdots \,|\, g_k(\theta_k)]$ by a constant so that its first entry is $1$, we obtain a Vandermonde matrix:\footnote{We use $r$ and $s$ to index the entries of a matrix, rather than $i$ and $j$. (We reserve $i$ to denote $\sqrt{-1}$.)}
\begin{align*}
\det(&[g_k(\theta_1) \,|\, \cdots \,|\, g_k(\theta_k)]) = \left(\prod_{s=1}^ke^{-\frac{k-1}{2}i\theta_s}\right)\det\!\left((e^{i\theta_s})^{r-1}\right)\!_{1 \le r,s \le k} \\
=& \prod_{s=1}^ke^{-\frac{k-1}{2}i\theta_s}\prod_{1 \le r < s \le k}(e^{i\theta_s} - e^{i\theta_r}) = \prod_{1 \le r < s \le k}e^{-\frac{i\theta_r}{2}}e^{-\frac{i\theta_s}{2}}(e^{i\theta_s} - e^{i\theta_r}) = \prod_{1 \le r < s \le k}2i\sin\!\left(\textstyle\frac{\theta_s - \theta_r}{2}\right)\!.
\end{align*}
We can check that $\pm(2i)^{-\lfloor k/2\rfloor}(2i)^{\binom{k}{2}} = 2^{\lfloor (k-1)^2/2\rfloor}$.
\end{pf}

\begin{prop}\label{fixed_point_properties}
The element $V_0\in\Gr_{k,n}(\mathbb{C})$ is represented by any $k\times n$ matrix of the form $[f_k(\theta + \frac{2\pi}{n}) \,|\, f_k(\theta + \frac{4\pi}{n}) \,|\, \cdots \,|\, f_k(\theta + 2\pi)]$, where $\theta\in\mathbb{R}$. Its Pl\"{u}cker coordinates are given by \rf{V_Plueckers}, and it is totally positive.
\end{prop}

\begin{pf}
Let $A_\theta$ denote the $k\times n$ matrix $[f_k(\theta + \frac{2\pi}{n}) \,|\, f_k(\theta + \frac{4\pi}{n}) \,|\, \cdots \,|\, f_k(\theta + 2\pi)]$. Then by \cref{vandermonde}, the $k\times k$ minor of $A_\theta$ located in columns $I = \{i_1 < \cdots < i_k\} \subseteq \{1, \cd, n\}$ equals
$$
\det(f_k(\theta + \textstyle\frac{2\pi i_1}{n}), \cd, f_k(\theta + \textstyle\frac{2\pi i_k}{n})) = 2^{\lfloor (k-1)^2/2\rfloor}\displaystyle\prod_{1 \le r < s \le k}\sin\!\left(\textstyle\frac{i_s - i_r}{n}\pi\right).
$$
Therefore $A_\theta$ represents a totally positive element $V$ of $\Gr_{k,n}(\mathbb{C})$ (independent of $\theta$), with Pl\"{u}cker coordinates given by \rf{V_Plueckers}. To see that $V$ equals $V_0$ (defined in \cref{defn_V(S)}), we take $\theta = -\frac{2\pi}{n}$, and apply row operations to $A_\theta$ as in \rf{exp_to_trig} to obtain a matrix with rows $(1, z, \cd, z^{n-1})$ for $z \in S_0$.
\end{pf}

We now establish a lemma which will imply that among the fixed points of $\sigma$, only $V_0$ can be totally nonnegative. The proof uses \cref{gantmakher_krein} (the result of Gantmakher and Krein).
\begin{lem}\label{argument_lemma}
Suppose that $V\in\Gr_{k,n}^{\ge 0}$, and $z\in\mathbb{C}^\times$ such that $(1, z, \cd, z^{n-1})\in V$. Then $|\!\arg(z)|\le\frac{k-1}{n-1}\pi$, where $\arg:\mathbb{C}^\times\to (-\pi, \pi]$ denotes the argument function.
\end{lem}

\begin{pf}
Let $\theta := |\!\arg(z)|$, and let $\epsilon\in\{1,-1\}$ be the sign of $\arg(z)$. Since $V$ has real Pl\"{u}cker coordinates, it is closed under complex conjugation. Hence for any $\phi\in\mathbb{R}$, the vector
$$
(\cos(\phi), |z|\cos(\theta+\phi), \cd, |z|^{n-1}\cos((n-1)\theta+\phi)) = \frac{e^{\epsilon i\phi}(1, z, \cd, z^{n-1}) + \overline{e^{\epsilon i\phi}(1, z, \cd, z^{n-1})}}{2}
$$
is in $V$. For $\phi<\frac{\pi}{2}$ sufficiently close to $\frac{\pi}{2}$, the vector on the left-hand side above changes sign exactly $\left\lceil\!\frac{(n-1)\theta}{\pi}\!\right\rceil$ times. Hence $\frac{(n-1)\theta}{\pi}\le k-1$ by \cref{gantmakher_krein}.
\end{pf}

\begin{rmk}
We note that the upper bound $\frac{k-1}{n-1}\pi$ of \cref{argument_lemma} is optimal. Indeed, given any $z\in\mathbb{C}^\times$ with $|\!\arg(z)|\le\frac{k-1}{n-1}\pi$, there exists $V\in\Gr_{k,n}^{\ge 0}$ with $(1, z, \cd, z^{n-1})\in V$. To see this, first note that by rescaling coordinate $j$ of $\mathbb{C}^n$ by $|z|^{-(j-1)}$ for $1 \le j \le n$, we may assume that $|z| = 1$. If $z=1$, we let $V$ be represented by the $k\times n$ matrix $(s^{r-1})_{1 \le r \le k, 1 \le s \le n}$. Otherwise, take $\rho := \frac{2|\!\arg(z)|}{k-1}$, and let $V$ be represented by the $k\times n$ matrix with columns $f_k(\theta)$, for $\theta = 0, \rho, 2\rho, \cd, (n-1)\rho$. Then $(1, z, \cd, z^{n-1})\in V$ by \rf{exp_to_trig}, and since $\rho\le\frac{2\pi}{n-1}$, we have $V\in\Gr_{k,n}^{\ge 0}$ by \cref{vandermonde}. (Note that if $\rho = \frac{2\pi}{n}$, then $V = V_0$.)
\end{rmk}

\begin{pf}[of \cref{unique_fixed_point}]
If $z$ is an $n$th root of $(-1)^{k-1}$, then $(1, z, \cd, z^{n-1})$ is an eigenvector of $\sigma\in\GL_n(\mathbb{C})$ with eigenvalue $z$. Thus $\sigma$ is diagonalizable with distinct eigenvalues, so a subspace $V\subseteq\mathbb{C}^n$ is $\sigma$-invariant if and only if it is spanned by eigenvectors of $\sigma$. This shows that the fixed points of $\sigma$ on $\Gr_{k,n}(\mathbb{C})$ are precisely $V_S$, for $k$-subsets $S\subseteq\mathbb{C}$ of $n$th roots of $(-1)^{k-1}$. These fixed points are all distinct by \cref{recover_S}. Now if $V_S$ is totally nonnegative, then \cref{argument_lemma} implies that $|\!\arg(z)| \le \frac{k-1}{n-1}\pi < \frac{k}{n}\pi$ for all $z\in S$ (if $k < n$), whence $S = S_0$ and $V_S = V_0$. We know from \rf{V_Plueckers} that $V_0$ is totally positive.
\end{pf}

\section{Quantum cohomology of the Grassmannian}\label{sec_quantum}

\noindent In this section we explain how the fixed points of the cyclic shift map appear in the theory of quantum cohomology. We discuss connections to the mirror-symmetric superpotential of the Grassmannian, and Schur polynomials evaluated at roots of unity.

\subsection{Peterson variety}
The {\itshape (small) quantum cohomology ring} $qH^*(\Gr_{k,n}(\mathbb{C}))$ of $\Gr_{k,n}(\mathbb{C})$ is a deformation of the cohomology ring $H^*(\Gr_{k,n}(\mathbb{C}))$ by an indeterminate $q$, defined as follows. Let $\mathcal{P}_{k,n}$ denote the set of partitions whose Young diagram fits inside the $k\times (n-k)$ rectangle, and for $\lambda\in\mathcal{P}_{k,n}$, let $\lambda^\vee\in\mathcal{P}_{k,n}$ denote the partition whose Young diagram (rotated by $180^\circ$) is the complement of the Young diagram of $\lambda$ inside the $k\times (n-k)$ rectangle. Then $H^*(\Gr_{k,n}(\mathbb{C}))$ has the {\itshape Schubert basis} $\{[X_\lambda] : \lambda\in\mathcal{P}_{k,n}\}$, where $[X_\lambda]$ is the cohomology class of the {\itshape Schubert variety} $X_\lambda$. By definition, $qH^*(\Gr_{k,n}(\mathbb{C}))$ has the basis $\{[X_\lambda] : \lambda\in\mathcal{P}_{k,n}\}$ over $\mathbb{Z}[q]$, with multiplication
$$
[X_\lambda]\star[X_\mu] := \sum_{d\ge 0}\sum_{\nu\in\mathcal{P}_{k,n}}\langle X_\lambda, X_\mu, X_{\nu^\vee}\rangle_d q^d [X_\nu].
$$
Here $\langle X_\lambda, X_\mu, X_{\nu^\vee}\rangle_d\in\mathbb{N}$ is a {\itshape Gromov--Witten invariant}, which is a generalized intersection number of the Schubert varieties $X_\lambda$, $X_\mu$, and $X_{\nu^\vee}$. If $d=0$ then the Gromov--Witten invariant is the usual intersection number given by the Littlewood--Richardson rule, and so the specialization at $q=0$ recovers the cup product in $H^*(\Gr_{k,n}(\mathbb{C}))$. {\itshape Quantum Schubert calculus} involves the study of these Gromov--Witten invariants; see \cite{bertram_97} for more details.

In unpublished work, Peterson defined a subvariety $\mathcal{Y}_{k,n}$ of $\GL_n(\mathbb{C})$ whose coordinate ring is isomorphic to $qH^*(\Gr_{k,n}(\mathbb{C}))$. This fact was proved by Rietsch \cite{rietsch01}, who characterized $\mathcal{Y}_{k,n}$ as follows.
\begin{defn}[Rietsch {\cite[Lemma 3.7 and (5.1)]{rietsch01}}]
For $z_1, \cd, z_k\in\mathbb{C}$, define the {\itshape Toeplitz matrix}
$$
u_{k,n}(z_1, \cd, z_k) := (e_{s-r}(z_1, \cd, z_k))_{1 \le r,s \le n}\in\GL_n(\mathbb{C}),
$$
where $e_j(x_1, \cd, x_k) := \sum_{1 \le i_1 < \cdots < i_j \le k}x_{i_1}\cdots x_{i_j}$ is the $j$th elementary symmetric polynomial for $j\ge 0$, and $e_j(x_1, \cd, x_k) := 0$ for $j < 0$. Then the {\itshape Peterson variety} is
$$
\mathcal{Y}_{k,n} := \{I_n\}\cup\{u_{k,n}(z_1, \cd, z_k) : z_1, \cd, z_k\in\mathbb{C}^\times\text{ are distinct and } z_1^n = \cdots = z_k^n\}\subseteq\GL_n(\mathbb{C}).
$$
Rietsch's isomorphism $qH^*(\Gr_{k,n}(\mathbb{C}))\to\mathbb{C}[\mathcal{Y}_{k,n}]$ identifies each element of the quantum cohomology ring with a function $\mathcal{Y}_{k,n}\to\mathbb{C}$. The indeterminate $q$ corresponds to the function
$$
I_n \mapsto 0, \quad u_{k,n}(z_1, \cd, z_k) \mapsto (-1)^{k-1}z_1^n,
$$
and the Schubert class $[X_\lambda]$ (for $\lambda\in\mathcal{P}_{k,n}$) corresponds to the function
$$
I_n \mapsto 0, \quad u_{k,n}(z_1, \cd, z_k) \mapsto s_\lambda(z_1, \cd, z_k).
$$
(Here $s_\lambda$ denotes the {\itshape Schur polynomial} indexed by $\lambda$, which we will study in \cref{Schur_subsection}.)
\end{defn}

There is a natural embedding $\gamma:\mathcal{Y}_{k,n}\hookrightarrow\Gr_{k,n}(\mathbb{C})$, whose image can be characterized using a deformation of the cyclic shift map $\sigma$.
\begin{defn}\label{deformed_defn}
For $t\in\mathbb{C}^\times$, define the {\itshape $t$-deformed (left) cyclic shift map} $\sigma_t\in\GL_n(\mathbb{C})$ by
$$
\sigma_t(v) := (v_2, v_3, \cd, v_n, (-1)^{k-1}tv_1) \quad \text{ for } v = (v_1, \cd, v_n)\in\mathbb{C}^n,
$$
so that $\sigma = \sigma_1$. Note that $\sigma_t$ induces an automorphism of $\Gr_{k,n}(\mathbb{C})$ of order $n$.
\end{defn}
We have an analogue of \cref{unique_fixed_point} for $\sigma_t$. We recall $V_S$, $S_0$, and $V_0$ from \cref{defn_V(S)}.
\begin{thm}\label{deformed_fixed_points}
For $t\in\mathbb{C}^\times$, the $t$-deformed cyclic shift map $\sigma_t$ on $\Gr_{k,n}(\mathbb{C})$ has exactly $\binom{n}{k}$ fixed points, namely, $V_S$ for any $k$-subset $S$ of distinct $n$th roots of $(-1)^{k-1}t$. Moreover, $V_S$ is totally nonnegative if and only if $t > 0$ and $S = t^{1/n}S_0$, in which case
\begin{align}\label{deformed_Plueckers}
\Delta_I(V_S) = t^{(i_1 + \cdots + i_s)/n}\prod_{1 \le r < s \le k}\sin\!\left(\textstyle\frac{i_s - i_r}{n}\pi\right) \quad \text{ for all $I = \{i_1 < \cdots < i_k\}\subseteq\{1, \cd, n\}$}.
\end{align}
\end{thm}

\begin{pf}
We only prove the second statement; the rest of the proof is similar to that of \cref{unique_fixed_point} at the end of \cref{sec_cyclic_shift}. If $t > 0$ and $S = t^{1/n}S_0$, then $V_S$ is obtained from $V_0$ by rescaling coordinate $j$ by $t^{j/n}$ for $1 \le j \le n$. Then \rf{deformed_Plueckers} follows from \rf{V_Plueckers}. Conversely, if $V_S$ is totally nonnegative, then it is closed under complex conjugation, whence so is $S$ by \cref{recover_S}. Therefore $t$ is real. Then \cref{argument_lemma} implies that $t > 0$ and $S = t^{1/n}S_0$.
\end{pf}

\begin{thm}\label{embedding_equality}
Define the embedding $\gamma:\mathcal{Y}_{k,n}\hookrightarrow\Gr_{k,n}(\mathbb{C}), g \mapsto \gamma(g)$ by
\begin{align}\label{embedding}
\Delta_I(\gamma(g)) = \det(g_{I^c,\{k+1, k+2, \cd, n\}}) \quad \text{ for all $k$-subsets }I\subseteq\{1, \cd, n\},
\end{align}
where $I^c$ denotes the complement of $I$. Then
\begin{align}\label{embedding_equation}
\gamma(u_{k,n}(z_1, \cd, z_k)) = V_{\{z_1, \cd, z_k\}} \quad \text{ for any distinct }z_1, \cd, z_k\in\mathbb{C}^\times \text{ with } z_1^n = \cdots = z_k^n.
\end{align}
Thus for any $t\in\mathbb{C}^\times$, by \cref{deformed_fixed_points} $\gamma$ sends the fiber of $q:\mathcal{Y}_{k,n}\to\mathbb{C}$ over $t$ to the set of fixed points of $\sigma_t$ in $\Gr_{k,n}(\mathbb{C})$.
\end{thm}

The proof uses the following lemma, which is classical. See \cite[Lemma 1.11(ii)]{karp_covectors} for a detailed discussion of this fact. Alternatively, one may prove \cref{embedding_equality} using the dual Jacobi--Trudi identity \cite[Corollary 7.16.2]{stanley_ec2} and \rf{Schur_evaluation}.
\begin{lem}\label{orthogonal}
Given $V\in\Gr_{k,n}(\mathbb{C})$, let $W\in\Gr_{n-k,n}(\mathbb{C})$ denote the subspace of $\mathbb{C}^n$ orthogonal to $V$ under the bilinear form $\langle v,w\rangle := v_1w_1 - v_2w_2 + \cdots + (-1)^{n-1}v_nw_n$ on $\mathbb{C}^n$. Then
$$
\Delta_I(V) = \Delta_{I^c}(W) \quad \text{ for all $k$-subsets }I\subseteq\{1, \cd, n\}.
$$
\end{lem}

\begin{pf}[of \cref{embedding_equality}]
We establish \rf{embedding_equation}. By \cref{orthogonal}, $\gamma(u_{k,n}(z_1, \cd, z_k))$ is the subspace of $\mathbb{C}^n$ orthogonal under $\langle\cdot,\cdot\rangle$ to the last $n-k$ columns of $u_{k,n}(z_1, \cd, z_k)$. Define the polynomial
$$
f(x) := \sum_{j=0}^ke_j(z_1, \cd, z_k)x^j = \prod_{j=1}^k(1+z_jx) \in \mathbb{C}[x].
$$
Note that for any $z\in\mathbb{C}^\times$ and $k+1 \le s \le n$, the pairing under $\langle\cdot,\cdot\rangle$ of $(1, z, \cd, z^{n-1})\in\mathbb{C}^n$ with column $s$ of $u_{k,n}(z_1, \cd, z_k)$ equals $(-z)^{s-1}f(-1/z)$, and hence vanishes when $z \in \{z_1, \cd, z_k\}$. That is, $(1, z, \cd, z^{n-1})$ lies in $\gamma(u_{k,n}(z_1, \cd, z_k))$ for $z\in\{z_1, \cd, z_k\}$.
\end{pf}

\begin{rmk}
We can verify from \cref{deformed_fixed_points} that as $t\to 0$, the fixed points of $\sigma_t$ converge to the single point $V\in\Gr_{k,n}(\mathbb{C})$ whose only nonzero Pl\"{u}cker coordinate is $\Delta_{\{1, \cd, k\}}$. (That is, any given neighborhood of $V$ contains all the fixed points of $\sigma_t$ for $|t|$ sufficiently small.) Therefore it is natural to call $V$ the unique fixed point of $\sigma_t$ at $t=0$, even though it does not make sense to set $t=0$ in \cref{deformed_defn}. This is consistent with \cref{embedding_equality}, since the fiber of $q : \mathcal{Y}_{k,n}\to\mathbb{C}$ over $0$ is $\{I_n\}$, and $\gamma(I_n) = V$.
\end{rmk}

\subsection{The mirror-symmetric superpotential}\label{sec_superpotential}

In studying mirror symmetry for partial flag varieties, Rietsch defined a rational function $\mathcal{F}_q$ on $\Gr_{k,n}(\mathbb{C})$, and proved the following:\footnote{Actually, $\mathcal{F}_q$ is a rational function on the Langlands dual of $\Gr_{k,n}(\mathbb{C})$, which we are identifying here with itself. The fact that $\Gr_{k,n}(\mathbb{C})$ is self-dual is a pleasant coincidence.}
\begin{thm}[Rietsch {\cite[Theorem 4.1]{rietsch08}}]\label{critical_points}
For $t\in\mathbb{C}^\times$, the set of critical points of $\mathcal{F}_q$ on $\Gr_{k,n}(\mathbb{C})$ at $q=t$ equals the fiber of $q:\mathcal{Y}_{k,n}\to\mathbb{C}$ over $t$, via the embedding $\gamma:\mathcal{Y}_{k,n}\hookrightarrow\Gr_{k,n}(\mathbb{C})$ from \rf{embedding}.
\end{thm}
In fact, Rietsch proved this result for general partial flag varieties $G/P$. Marsh and Rietsch further studied the case $G/P = \Gr_{k,n}(\mathbb{C})$, and gave the following simple formula for $\mathcal{F}_q$ (which they call the {\itshape superpotential}) \cite[(6.4)]{marsh_rietsch}:
\begin{align}\label{superpotential_formula}
\mathcal{F}_q = \left(\sum_{1 \le i \le n, \;\, i\neq n-k}\frac{\Delta_{\{i+1, i+2, \cd, i+k-1, i+k+1\}}}{\Delta_{\{i+1, i+2, \cd, i+k\}}}\right) + q\frac{\Delta_{\{n-k+1, n-k+2, \cd, n-1, 1\}}}{\Delta_{\{n-k+1, n-k+2, \cd, n\}}}.
\end{align}
Here the indices are taken modulo $n$. From \cref{embedding_equality} and \cref{critical_points}, we find:
\begin{cor}\label{critical_equals_fixed}
For $t\in\mathbb{C}^\times$, the critical points of $\mathcal{F}_q$ on $\Gr_{k,n}(\mathbb{C})$ at $q=t$ are precisely the fixed points of the $t$-deformed cyclic shift map $\sigma_t$.
\end{cor}
It may also be possible to establish this correspondence (but only for a subset of the critical points) via {\itshape birational rowmotion}; see \cref{sec_rowmotion}.

\begin{eg}\label{superpotential_eg}
Let $k := 2$, $n := 4$. Then representing a generic point $V\in\Gr_{2,4}(\mathbb{C})$ by the matrix $\scalebox{0.9}{$\begin{bmatrix}
1 & 0 & -b & -\frac{bc}{a} - d \\[3pt]
0 & 1 & a & c
\end{bmatrix}$}$ (see \rf{laurent_superpotential} for motivation behind this choice), we obtain
$$
\qquad\mathcal{F}_q = \frac{\Delta_{\{1,3\}}}{\Delta_{\{1,2\}}} + \frac{\Delta_{\{2,4\}}}{\Delta_{\{2,3\}}} + q\frac{\Delta_{\{1,3\}}}{\Delta_{\{3,4\}}} + \frac{\Delta_{\{2,4\}}}{\Delta_{\{1,4\}}} = a + \frac{c}{a} + \frac{d}{b} + q\frac{1}{d} + \frac{b}{a} + \frac{d}{c}.
$$
Let us set $q=1$. We can check that $\mathcal{F}_q$ has critical points $V$ corresponding to
$$
(a,b,c,d) = (\sqrt{2}, 1, 1, \textstyle\frac{1}{\sqrt{2}}),\; (-\sqrt{2}, 1, 1, -\textstyle\frac{1}{\sqrt{2}}),\; (\sqrt{2}i, -1, -1, -\textstyle\frac{i}{\sqrt{2}}),\; (-\sqrt{2}i, -1, -1, \textstyle\frac{i}{\sqrt{2}}).
$$
There are two additional critical points, occurring when $\Delta_{\{1,3\}} = \Delta_{\{2,4\}} = 0$ (which we do not see in our choice of coordinates), giving $6$ in total. By \cref{critical_equals_fixed}, these are precisely the $6$ fixed points of the cyclic shift map $\sigma$. The point $(a,b,c,d) = (\sqrt{2}, 1, 1, \frac{1}{\sqrt{2}})$ corresponds to the unique totally positive fixed point $V_0$ of $\sigma$ from \cref{intro_eg}.
\end{eg}

\subsection{Schur polynomials evaluated at roots of unity}\label{Schur_subsection}

The {\itshape Schur polynomial} of a partition $\lambda = (\lambda_1, \cd, \lambda_k)$ with at most $k$ parts is defined as
\begin{align}\label{Schur_defn}
s_\lambda(x_1, \cd, x_k) := \frac{\det(x_r^{\lambda_{k+1-s}+s-1})_{1 \le r,s \le k}}{\det(x_r^{s-1})_{1 \le r,s \le k}}.
\end{align}
In particular, for any distinct $z_1, \cd, z_k\in\mathbb{C}$, we have
\begin{align}\label{Schur_evaluation}
s_\lambda(z_1, \cd, z_k) = \frac{\Delta_{\{\lambda_k + 1, \lambda_{k-1} + 2, \cd, \lambda_1 + k\}}(V_{\{z_1, \cd, z_k\}})}{\Delta_{\{1, \cd, k\}}(V_{\{z_1, \cd, z_k\}})} \quad \text{ for all }\lambda\in\mathcal{P}_{k,n}.
\end{align}
Then reinterpreting the second statement of \cref{deformed_fixed_points} in terms of Schur polynomials, we recover the following result of Rietsch, which she proved using orthogonality relations for Schur polynomials.
\begin{cor}[Rietsch {\cite[Theorem 8.4]{rietsch01}}]\label{Schur_corollary}
Let $0 \le k < n$, $t\in\mathbb{C}^\times$, and $z_1, \cd, z_k\in\mathbb{C}$ be distinct $n$th roots of $(-1)^{k-1}t$. Then $s_\lambda(z_1, \cd, z_k)\ge 0$ for all $\lambda\in\mathcal{P}_{k,n}$ if and only if $t > 0$ and $\{z_1, \cd, z_k\} = t^{1/n}S_0$. In this case,
\begin{align}\label{Schur_Plueckers}
s_\lambda(z_1, \cd, z_k) = t^{|\lambda|/n}\prod_{1 \le r < s \le k}\frac{\sin\left(\textstyle\frac{\lambda_r-\lambda_s+s-r}{n}\pi\right)}{\sin\left(\textstyle\frac{s-r}{n}\pi\right)} \quad \text{ for all }\lambda\in\mathcal{P}_{k,n}.
\end{align}
\end{cor}
We remark that we can instead use Rietsch's result to deduce the second statement of \cref{deformed_fixed_points}, as an alternative to using \cref{argument_lemma}. We also mention that in their study of symmetric group characters, Orellana and Zabrocki \cite{orellana_zabrocki} consider evaluations of symmetric polynomials at certain other $k$-multisubsets of roots of unity (namely, those which are the eigenvalues of a $k\times k$ permutation matrix).

\begin{rmk}
Rietsch established \cref{Schur_corollary} in order to describe the {\itshape totally nonnegative part $\mathcal{Y}_{k,n}^{\ge 0}$} of $\mathcal{Y}_{k,n}$, which is defined as the subset of $\mathcal{Y}_{k,n}$ of matrices whose minors are all nonnegative. Now recall the embedding $\mathcal{Y}_{k,n}\hookrightarrow\Gr_{k,n}(\mathbb{C})$ from \rf{embedding}. It is clear that if $g\in\mathcal{Y}_{k,n}^{\ge 0}$, then $\gamma(g)\in\Gr_{k,n}^{\ge 0}$. Rietsch showed that the converse holds, using the combinatorics of Schur polynomials. Therefore \cref{embedding_equality} and \cref{deformed_fixed_points} imply that
$$
\mathcal{Y}_{k,n}^{\ge 0} = \gamma^{-1}(\Gr_{k,n}^{\ge 0}) = \{I_n\} \cup \{u_{k,n}(t^{1/n}S_0) : t > 0\},
$$
recovering \cite[Theorem 8.4(2)]{rietsch01}. Rietsch also gives an explicit factorization of the element $u_{k,n}(t^{1/n}S_0)$ into elementary matrices \cite[Proposition 9.3]{rietsch01}.
\end{rmk}

\begin{rmk}
Rietsch \cite[Proposition 11.1]{rietsch01} proved the following intriguing inequality: for any distinct $z_1, \cd, z_k\in\mathbb{C}$ such that $z_1^n = \cdots = z_k^n$ and $|z_1| = 1$, we have
$$
|s_\lambda(z_1, \cd, z_k)| \le s_\lambda(S_0) = \prod_{1 \le r < s \le k}\frac{\sin\left(\textstyle\frac{\lambda_r-\lambda_s+s-r}{n}\pi\right)}{\sin\left(\textstyle\frac{s-r}{n}\pi\right)}\quad\text{ for all }\lambda\in\mathcal{P}_{k,n}.
$$
It would be interesting to explain how these inequalities are related to the discrete isoperimetric inequalities of Nudel'man discussed in \cref{sec_isoperimetric_inequalities}.
\end{rmk}

\section{Related work}\label{sec_connections}

\subsection{Cyclic and bicyclic polytopes}\label{sec_cyclic_polytopes}

Recall the curves $f_k$ defined in \rf{f_defn}. For odd $k$, the curve in $\mathbb{R}^{k-1}$ formed from $f_k$ by deleting the first component (or by regarding $f_k$ as a curve in $\mathbb{P}^{k-1}$) is the {\itshape trigonometric moment curve}. Carath\'{e}odory \cite{caratheodory} used such curves, along with the {\itshape (power) moment curves} $t\mapsto (t, t^2, \cd, t^{k-1})$, to define cyclic polytopes of even dimension. Namely, a polytope is {\itshape cyclic} if it has the same face lattice as a polytope whose vertices lie on such a curve. Cyclic polytopes have the maximum number of faces in each dimension among all simplicial $(k-2)$-spheres with a fixed number of vertices; this is the celebrated upper bound theorem of McMullen \cite{mcmullen} and Stanley \cite{stanley_75}.

For even $k$, the curve $f_k$ is the {\itshape symmetric moment curve}, first studied by Nudel'man \cite{nudel'man} in order to resolve an isoperimetric problem (see \cref{sec_isoperimetric_inequalities}). Motivated by establishing an upper bound theorem for centrally symmetric polytopes (which is an open problem, even for $4$-dimensional polytopes), Barvinok and Novik \cite{barvinok_novik} employed symmetric moment curves to define a centrally symmetric analogue of cyclic polytopes, called {\itshape bicyclic polytopes}. They used bicyclic polytopes to give lower bounds on the numbers of faces of centrally symmetric polytopes, which they and Lee \cite{barvinok_lee_novik_3} improved with polytopes constructed from variants of the symmetric moment curve. We refer to the survey of Novik \cite{novik19} for a detailed discussion of this topic.

We remark that before Barvinok and Novik's work \cite{barvinok_novik}, there was not a consensus about what the centrally symmetric analogue of the cyclic polytope should be. We view \cref{unique_fixed_point} as further evidence that the symmetric moment curve belongs in the same family as the trigonometric moment curve.

\subsection{Isoperimetric inequalities for convex curves}\label{sec_isoperimetric_inequalities} A continuous curve $g : S\to\mathbb{R}^k$ defined on an interval $S\subseteq\mathbb{R}$ is called {\itshape convex on $\mathbb{R}^k$} if the determinants
$$
\det\!\left(\begin{bmatrix}1 & 1 & \cdots & 1 \\ g(t_0) & g(t_1) & \cdots & g(t_k)\end{bmatrix}\right) \qquad (t_0 < t_1 < \cdots < t_k \text{ in } S)
$$
are either all nonnegative or all nonpositive, and not all zero.\footnote{If we also require that the determinants are nonzero, we obtain the definition of a {\itshape strictly convex} curve. Such curves have been widely studied under the various names {\itshape curves of order $k$} \cite{juel}, {\itshape monotone curves} \cite{hjelmslev}, {\itshape strictly comonotone curves} \cite{motzkin}, {\itshape alternating curves} \cite{bjorner_las_vergnas_sturmfels_white_ziegler}, and {\itshape hyperconvex curves} \cite{labourie}.} The convex hull of such a curve is the continuous analogue of a cyclic polytope. Indeed, recall that a polytope is {\itshape cyclic} if its face lattice is the same as a polytope whose vertices lie on the curve $t\mapsto (t, t^2, \cd, t^k)\in\mathbb{R}^k$; we call it {\itshape alternating} if in addition every subset of its vertices generates a cyclic polytope.\footnote{Shemer \cite[Theorem 2.12]{shemer_82} showed that in even dimension, every cyclic polytope is alternating. In odd dimension, see \cite[pp.\ 396-397]{bjorner_las_vergnas_sturmfels_white_ziegler} for an example of a cyclic polytope which is not alternating.} Sturmfels \cite{sturmfels_88} showed that $x_1, \cd, x_n\in\mathbb{R}^k$ are the vertices of a $k$-dimensional alternating polytope if and only if the matrix \scalebox{0.7}{$\begin{bmatrix}1 & \cdots & 1 \\ x_1 & \cdots & x_n\end{bmatrix}$} represents a totally positive element of $\Gr_{k+1,n}(\mathbb{R})$. In this case, the polygonal path through $x_1, \cd, x_n$ is convex on $\mathbb{R}^k$, so every alternating polytope is the convex hull of a convex curve. Every convex curve has the remarkable property \cite[Theorem II.3.1]{karlin_studden}
that its {\itshape Carath\'{e}odory number} is $\lfloor k/2\rfloor+1$, i.e.\ any point in its convex hull is the convex combination of at most $\lfloor k/2\rfloor+1$ points on the curve. In general, the Carath\'{e}odory number of a curve in $\mathbb{R}^k$ is at most $k$, by a result of Fenchel \cite{fenchel29}.

Let $L$ denote the length of a curve convex on $\mathbb{R}^k$, and $V$ the volume of its convex hull. The {\itshape isoperimetric problem for convex curves}, posed by Schoenberg \cite{schoenberg_54}, is to find (given $k$) an upper bound for $V$ in terms of $L$. There are three cases, depending on the parity of $k$ and whether the curve is closed or not. (There is no fourth case, since closed convex curves only exist in even dimension.) Note the similarities of the extremal curves to $f_{k+1}$ or $f_k$.\vspace*{4pt}
\begin{center}
\small\setlength{\tabcolsep}{5pt}
\begin{tabular}{|m{1.08in}|>{\vspace*{3pt}}m{1.42in}<{\vspace*{3pt}}|m{3.54in}|}
\hline
{\bfseries\parbox{1.08in}{\centering case}} & {\bfseries\parbox{1.42in}{\centering isoperimetric \\ inequality}} & {\bfseries\parbox{3.54in}{\centering extremal curve, for $\theta\in [0,2\pi]$}} \\ \hline
\parbox{1.08in}{$k$ even; closed \\ curves \cite{schoenberg_54}} & $V\le\displaystyle\frac{L^k}{(\pi k)^{\frac{k}{2}}k!(\frac{k}{2})!}$ & $\left(\!\cos(\theta), \sin(\theta), \frac{\cos(2\theta)}{2}, \frac{\sin(2\theta)}{2}, \cd, \frac{\cos(\frac{k}{2}\theta)}{\frac{k}{2}}, \frac{\sin(\frac{k}{2}\theta)}{\frac{k}{2}}\!\right)$ \\ \hline
\parbox{1.09in}{\begin{flushleft}$k$ odd \cite[Theorem III.8.6]{krein_nudel'man}\end{flushleft}} & $V\le\displaystyle\frac{L^k}{\pi^{\frac{k-1}{2}}k^{\frac{k}{2}}k!(\frac{k-1}{2})!}$ & \scalebox{0.97}{$\left(\!\frac{\theta}{\sqrt{2}}, \cos(\theta), \sin(\theta), \frac{\cos(2\theta)}{2}, \frac{\sin(2\theta)}{2}, \cd, \frac{\cos(\frac{k-1}{2}\theta)}{\frac{k-1}{2}}, \frac{\sin(\frac{k-1}{2}\theta)}{\frac{k-1}{2}}\!\right)$} \\ \hline
\parbox{1.08in}{$k$ even \cite{nudel'man}} & $V\le\displaystyle\frac{L^k}{(\frac{\pi k}{2})^{\frac{k}{2}}k!(k-1)!!}$ & \scalebox{0.94}{$\left(\!\cos\left(\frac{1}{2}\theta\right), \sin\left(\frac{1}{2}\theta\right), \frac{\cos(\frac{3}{2}\theta)}{3}, \frac{\sin(\frac{3}{2}\theta)}{3}, \cd, \frac{\cos(\frac{k-1}{2}\theta)}{k-1}, \frac{\sin(\frac{k-1}{2}\theta)}{k-1}\!\right)$} \\ \hline
\end{tabular}
\end{center}
\vspace*{4pt}In the first two cases above, the extremal curve is unique (modulo affine isometries of $\mathbb{R}^k$). It is not known if uniqueness holds in the third case.

There is a discrete version of this problem, in which we fix $n\ge k$ and consider curves convex on $\mathbb{R}^k$ that are polygonal paths with $n$ segments. This is the isoperimetric problem for alternating polytopes. It was proposed by Krein and Nudel'man (see the discussion after \cite[Theorem III.8.6]{krein_nudel'man}) and solved by Nudel'man
\cite{nudel'man}:\vspace*{4pt}
\begin{center}
\small\setlength{\tabcolsep}{5pt}
\begin{tabular}{|m{0.85in}|>{\vspace*{3pt}}m{1.99in}<{\vspace*{3pt}}|m{3.18in}|}
\hline
{\bfseries\parbox{0.85in}{\centering case}} & {\bfseries\parbox{1.99in}{\centering isoperimetric \\ inequality}} & {\bfseries\parbox{3.18in}{\centering vertices of extremal piecewise linear curve, for $\theta = 0, \frac{2\pi}{n}, \frac{4\pi}{n}, \cd, 2\pi$}} \\ \hline
\parbox{0.85in}{\vspace*{3pt} $k$ even, \\ $n$ segments; \\ closed curves\vspace*{3pt}} & $V\le\displaystyle\frac{L^k}{k^{\frac{k}{2}}k!\prod_{j=1}^{\frac{k}{2}}n\tan(\frac{j\pi}{n})}$ & $\left(\!\frac{\cos(\theta)}{\sin(\frac{\pi}{n})}, \frac{\sin(\theta)}{\sin(\frac{\pi}{n})}, \frac{\cos(2\theta)}{\sin(\frac{2\pi}{n})}, \frac{\sin(2\theta)}{\sin(\frac{2\pi}{n})}, \cd, \frac{\cos(\frac{k}{2}\theta)}{\sin(\frac{k\pi}{2n})}, \frac{\sin(\frac{k}{2}\theta)}{\sin(\frac{k\pi}{2n})}\!\right)$ \\ \hline
\parbox{0.85in}{$k$ odd, \\ $n$ segments} & $V\le\displaystyle\frac{L^k}{k^{\frac{k}{2}}k!\prod_{j=1}^{\frac{k-1}{2}}n\tan(\frac{j\pi}{n})}$ & \scalebox{0.82}{$\left(\!\frac{n\theta}{\sqrt{2}\pi}, \frac{\cos(\theta)}{\sin(\frac{\pi}{n})}, \frac{\sin(\theta)}{\sin(\frac{\pi}{n})}, \frac{\cos(2\theta)}{\sin(\frac{2\pi}{n})}, \frac{\sin(2\theta)}{\sin(\frac{2\pi}{n})}, \cd, \frac{\cos(\frac{k-1}{2}\theta)}{\sin(\frac{(k-1)\pi}{2n})}, \frac{\sin(\frac{k-1}{2}\theta)}{\sin(\frac{(k-1)\pi}{2n})}\!\right)$} \\ \hline
\parbox{0.86in}{$k$ even, \\ $n$ segments\footnotemark} & $V\le\displaystyle\frac{L^k}{k^{\frac{k}{2}}k!\prod_{j=1}^{\frac{k}{2}}n\tan(\frac{(2j-1)\pi}{2n})}$ & \scalebox{0.86}{$\left(\!\frac{\cos(\frac{1}{2}\theta)}{\sin(\frac{\pi}{2n})}, \frac{\sin(\frac{1}{2}\theta)}{\sin(\frac{\pi}{2n})}, \frac{\cos(\frac{3}{2}\theta)}{\sin(\frac{3\pi}{2n})}, \frac{\sin(\frac{3}{2}\theta)}{\sin(\frac{3\pi}{2n})}, \cd, \frac{\cos(\frac{k-1}{2}\theta)}{\sin(\frac{(k-1)\pi}{2n})}, \frac{\sin(\frac{k-1}{2}\theta)}{\sin(\frac{(k-1)\pi}{2n})}\!\right)$} \\ \hline
\end{tabular}\footnotetext{We have corrected a typo in \cite[(2)]{nudel'man}.}
\end{center}
\vspace*{4pt}In all three cases, the extremal curve is unique. By taking $n\to\infty$, we recover the isoperimetric inequalities in the continuous case, but cannot conclude that the resulting extremal curve is unique. We note that Nudel'man's constructions are similar to the definition of $V_0$, but it is difficult to make this observation more precise. This is because a real matrix representing $V_0$ is unique modulo linear automorphisms of $\mathbb{R}^k$, while the extremal curves above are unique modulo affine isometries of $\mathbb{R}^k$.

We again make a remark that before Nudel'man's work \cite{nudel'man}, it was not known what the extremal curve for even $k$ would be. Krein and Nudel'man \cite{krein_nudel'man} speculated it might be half of the extremal curve in the closed case, which is true for $k=2$ but false in general.

\subsection{Separable quantum states}\label{sec_separable_states}

Recall that the trigonometric moment curve in $\mathbb{R}^{2d}$ is formed from $f_{2d+1}$ (see \rf{f_defn}) by deleting the first component. Let $C_{2d}\subseteq\mathbb{R}^{2d}$ denote its convex hull. Kye \cite{kye13} found a surprising appearance of $C_{2d}$ in quantum mechanics.

We regard the space of $d\times d$ complex matrices as the C*-algebra of bounded operators $\mathcal{B}(\mathbb{C}^d)$ on $\mathbb{C}^d$. A {\itshape state} on $\mathcal{B}(\mathbb{C}^d)$ is a linear functional $\phi:\mathcal{B}(\mathbb{C}^d)\to\mathbb{C}$ which takes the identity to $1$ and satisfies $\phi(A)\ge 0$ if $A$ is positive semidefinite. Equivalently, there exists a positive semidefinite $B\in\mathcal{B}(\mathbb{C}^d)$ with $\tr(B) = 1$ such that $\phi(A) = \tr(AB)$ for all $A\in\mathcal{B}(\mathbb{C}^d)$ (this is a special case of the Riesz representation theorem). Note that for any positive semidefinite $B\in\mathcal{B}(\mathbb{C}^d)$ and $B'\in\mathcal{B}(\mathbb{C}^{d'})$, we obtain a state on $\mathcal{B}(\mathbb{C}^d)\otimes\mathcal{B}(\mathbb{C}^{d'})\cong\mathcal{B}(\mathbb{C}^d\otimes\mathbb{C}^{d'})$ given by $A\otimes A'\mapsto \tr(AB)\tr(A'B')$. A convex combination of such states on $\mathcal{B}(\mathbb{C}^d)\otimes\mathcal{B}(\mathbb{C}^{d'})$ is called {\itshape separable}; non-separable states are called {\itshape entangled}. The space of separable states $\mathbb{S}_{d,d'}$ on $\mathcal{B}(\mathbb{C}^d)\otimes\mathcal{B}(\mathbb{C}^{d'})$ is a real, compact, convex set of dimension $(dd')^2-1$. In this context, Kye showed the following:
\begin{thm}[Kye \cite{kye13}]
$C_{2d}$ is a face of $\mathbb{S}_{d,2}$.
\end{thm}

Kye also recovered a description of Carath\'e{o}dory \cite[Section 12]{caratheodory} of the faces of $C_{2d}$: they are precisely the simplices generated by some $d$ or fewer distinct points on the trigonometric moment curve, and are all exposed.\footnote{A {\itshape face} of a convex set $C$ is a convex subset $F$ such that any line segment in $C$ whose interior intersects $F$ is contained entirely in $F$, and $F$ is {\itshape exposed} if it is the intersection of $C$ with an affine hyperplane.} This essentially follows from the fact that the Carath\'e{o}dory number of the trigonometric moment curve is $d+1$; see \cite[Corollary 5.4]{sanyal_sottile_sturmfels11}.

We remark that the facial structure of the convex hull of the symmetric moment curve $f_{2d}$ in $\mathbb{R}^{2d}$ has been studied by Smilansky (in the case $d=2$) \cite{smilansky85}, Barvinok, Lee, and Novik \cite{barvinok_novik, barvinok_lee_novik_1}, and Vinzant \cite{vinzant11}.

\subsection{Arrangements of equal minors}\label{arrangements_of_equal_minors}
Farber and Postnikov \cite{farber_postnikov} studied the possible arrangements of equal and unequal Pl\"{u}cker coordinates (or ``minors'') among totally positive elements of $\Gr_{k,n}(\mathbb{C})$. By \rf{V_Plueckers}, $V_0$ has many pairs of equal Pl\"{u}cker coordinates: if $I$ and $J$ are $k$-subsets of $\{1, \cd, n\}$ which are cyclic shifts or reflections of each other modulo $n$, then $\Delta_I(V_0) = \Delta_J(V_0)$. (This does not hold for all $V\in\Gr_{k,n}(\mathbb{C})$ fixed by $\sigma$: we have $\sigma(V) = V$ if and only if there exists $\zeta\in\mathbb{C}$ with $\zeta^n = 1$ such that $\Delta_{\{i_1+j, \cd, i_k+j\}}(V) = \zeta^j\Delta_{\{i_1, \cd, i_k\}}(V)$ for all $j\in\mathbb{Z}$ and distinct $i_1, \cd, i_k$ modulo $n$.) We do not know if the converse is true, i.e.\ if $\Delta_I(V_0) = \Delta_J(V_0)$ implies that $I$ and $J$ are related by dihedral symmetry.

Farber and Postnikov were especially interested in the minimum and maximum Pl\"{u}cker coordinates of totally positive elements.
\begin{prob}
What are the minimum and maximum Pl\"{u}cker coordinates of $V_0$ (after we have rescaled its Pl\"{u}cker coordinates to be strictly positive)?
\end{prob}
We expect that the minimum Pl\"{u}cker coordinates are indexed by $\{1, \cd, k\}$ and its cyclic shifts. In general, $\Delta_I(V_0)$ should measure `how spread out' are the elements of $I$ modulo $n$. We observe that we can interpret $\Delta_I(V_0)$ as the volume of a simplex generated by vertices of an alternating polytope.

\subsection{The twist map and its fixed points}\label{sec_twist}

The Grassmannian {\itshape twist map} was introduced by Marsh and Scott \cite{marsh_scott}, as an analogue of the twist map on unipotent matrices defined by Berenstein, Fomin, and Zelevinsky \cite[(3.1.1)]{berenstein_fomin_zelevinsky96}. It is important in the study of the cluster-algebraic structure of the Grassmannian, connecting the {\itshape $\mathcal{A}$-cluster structure} and {\itshape $\mathcal{X}$-cluster structure} of the homogeneous coordinate ring of $\Gr_{k,n}(\mathbb{C})$, and is closely related to the {\itshape Donaldson--Thomas transformation} considered by Weng \cite{weng}. The twist map was generalized from $\Gr_{k,n}(\mathbb{C})$ to all its {\itshape positroid subvarieties} by Muller and Speyer \cite{muller_speyer}.

The {\itshape (right) twist map} $\tau$ is an automorphism of $\Pi_{k,n}^\circ$, where $\Pi_{k,n}^\circ$ is the subset of $\Gr_{k,n}(\mathbb{C})$ in which the $n$ Pl\"{u}cker coordinates indexed by cyclic intervals of $\{1, \cd, n\}$ are nonzero.\footnote{$\Pi_{k,n}^\circ$ is the top-dimensional stratum of the {\itshape positroid stratification} of $\Gr_{k,n}(\mathbb{C})$ studied by Knutson, Lam, and Speyer \cite{knutson_lam_speyer}. This stratification is the projection to $\Gr_{k,n}(\mathbb{C})$ of the stratification of the complete flag variety into intersections of opposite Schubert cells, which was first considered by Lusztig \cite{lusztig}. Rietsch \cite{rietsch99} proved his conjecture that it induces a cell decomposition of $\Gr_{k,n}^{\ge 0}$. Muller and Speyer extended the definition of the twist map to all of $\Gr_{k,n}(\mathbb{C})$, by defining it individually on each positroid stratum.} We use Muller and Speyer's definition of $\tau$, which differs from Marsh and Scott's by rescaling coordinates. Take the standard bilinear form $\langle v,w\rangle := v_1w_1 + \cdots + v_kw_k$ on $\mathbb{C}^k$, and let $V\in\Pi_{k,n}^\circ$ be represented by the $k\times n$ matrix $A$. Then $\tau(V)$ is represented by the $k\times n$ matrix whose $j$th column is orthogonal to columns $j+1, \cd, j+k-1$ (modulo $n$) of $A$, and pairs to $1$ with column $j$ of $A$. For example, $\tau$ sends $\scalebox{0.9}{$\begin{bmatrix}1 & 1 & 0 & -4 \\[2pt] 0 & 2 & 1 & 3\end{bmatrix}$}$ to $\scalebox{0.9}{$\begin{bmatrix}1 & 1 & \frac{3}{4} & 0 \\[2pt] -\frac{1}{2} & 0 & 1 & \frac{1}{3}\end{bmatrix}$}$. The inverse of the right twist $\tau$ is the analogously-defined {\itshape left twist}.

We remark that Marsh and Rietsch employed the twist map to study the superpotential $\mathcal{F}_q$ \cite[Section 12]{marsh_rietsch}. This uses an elegant expression of Marsh and Scott \cite[Proposition 3.5]{marsh_scott} for the twist of certain Pl\"{u}cker coordinates; these coordinates happen to include those appearing in \rf{superpotential_formula}. The twist map $\tau$ is also closely related to the cyclic shift map $\sigma$, after we quotient by the {\itshape torus action} on $\Gr_{k,n}(\mathbb{C})$:
\begin{thm}[Marsh and Scott {\cite[Corollary 4.2]{marsh_scott}}]\label{periodicity}
Let $(\mathbb{C}^\times)^n$ act on $\Gr_{k,n}(\mathbb{C})$ by rescaling coordinates, i.e.\ rescaling columns of $k\times n$ matrices. Then $\tau^2 = \sigma^k$ on $\Pi_{k,n}^\circ/(\mathbb{C}^\times)^n$.
\end{thm}
In fact, one of our motivations for studying $\sigma$ was to try to better understand $\tau$.

\cref{periodicity} implies that $\tau$ is periodic of order $2n$ on $\Pi_{k,n}^\circ/(\mathbb{C}^\times)^n$, and that every fixed point of $\sigma$ in $\Pi_{k,n}^\circ/(\mathbb{C}^\times)^n$ is fixed by $\tau^2$. Here we observe that many of these points are also fixed by $\tau$. We work more generally over $\Pi_{k,n}^\circ$, since our method does not identify any fixed points of $\tau$ in $\Pi_{k,n}^\circ/(\mathbb{C}^\times)^n$ which do not come from fixed points in $\Pi_{k,n}^\circ$.
\begin{prop}
Let $S = \{z_1, \cd, z_k\}\subseteq\mathbb{C}^\times$ such that $z_1^n = \cdots = z_k^n$ and $S$ is closed under inversion (i.e.\ $S = \{z_1^{-1}, \cd, z_k^{-1}\}$). Then the element $V_S$ of $\Gr_{k,n}(\mathbb{C})$ (defined in \cref{defn_V(S)}) is fixed by the twist map $\tau$. Moreover, $V_S$ is totally nonnegative if and only if $V_S = V_0$.
\end{prop}

\begin{pf}
We represent $V_S$ by the $k\times n$ matrix $A := (z_r^{s-1})_{1 \le r \le k, 1 \le s \le n}$. For $1 \le j \le n$, define the $k\times k$ matrix $B_j := (z_r^{j+s-2})_{1 \le r,s \le k}$. Note that column $1$ of $B_j$ is column $j$ of $A$, and since $z_1^n = \cdots = z_k^n$, columns $2, \cd, k$ of $B_j$ are columns $j+1, \cd, j+k-1$ (modulo $n$) of $A$, up to rescaling. Therefore $\tau(V_S)$ is represented by the $k\times n$ matrix $A'$ whose $j$th column equals the first row of $B_j^{-1}$. Applying the classical adjoint description of the matrix inverse and Vandermonde's determinantal identity, we find that the $(p,j)$-entry of $A'$ equals
$$
(-1)^{p-1}\frac{\det(z_r^{j+s-2})_{\substack{r = 1, \cd, p-1, p+1, \cd, k \\ s = 2, \cd, k\phantom{,p-1, \cd, p+1}}}}{\det(z_r^{j+s-2})_{1 \le r,s \le k}} = (-1)^{p-1}\frac{(z_1\cdots z_k)^jz_p^{-j}\displaystyle\prod_{r<s,\;\, r,s\neq p}(z_s - z_r)}{(z_1\cdots z_k)^{j-1}\displaystyle\prod_{r<s}(z_s - z_r)} = \frac{(z_1\cdots z_k)z_p^{-j}}{\displaystyle\prod_{s\neq p}(z_s - z_p)}.
$$
(This calculation also shows that $V_S$ is in $\Pi_{k,n}^\circ$.) Since $S$ is closed under inversion, $A'$ equals $A$ up to rescaling and permuting the rows. The last statement follows from \cref{deformed_fixed_points}.
\end{pf}

\begin{prob}
Find all fixed points of the twist map $\tau$ on $\Pi_{k,n}^\circ$, and more generally, on any positroid stratum of $\Gr_{k,n}(\mathbb{C})$. Is $V_0$ the unique totally positive fixed point?
\end{prob}

\subsection{Promotion on rectangular tableaux}\label{sec_promotion}

Fix $k\le n$, and let $T$ be a semistandard tableau on $\{1, \cd, n\}$ of rectangular shape with $n-k$ rows (and any number of columns). The {\itshape promotion} $\pr(T)$ of $T$, first studied by Sch\"{u}tzenberger \cite{schutzenberger72}, is the semistandard tableau of the same shape formed by cyclically shifting the entries (sending $i$ to $i-1$ for $2 \le i \le n$, and $1$ to $n$), and semistandardizing via {\itshape jeu-de-taquin slides}. For example (with $k=3$, $n=5$),
$$
T = \begin{ytableau}1 & 1 & 2 & 3 \\ 2 & 3 & 4 & 5\end{ytableau} \;\;\rightsquigarrow\;\; \begin{ytableau}5 & 5 & 1 & 2 \\ 1 & 2 & 3 & 4\end{ytableau} \;\;\rightsquigarrow\;\; \begin{ytableau}1 & 1 & 2 & 4 \\ 2 & 3 & 5 & 5\end{ytableau}\, = \pr(T).
$$
In constructing an affine geometric crystal on the Grassmannian, Frieden \cite[Theorem 5.4]{frieden19} showed that the deformed cyclic shift map $\sigma_t$ on $\Gr_{k,n}(\mathbb{C})$ {\itshape tropicalizes} to promotion $\pr$, in a very precise sense.\footnote{Our conventions differ from Frieden's: our ground set $\{1, \cd, n\}$ has the opposite order to his.} The superpotential $\mathcal{F}_q$ also appears, as the {\itshape decoration} of the crystal.

We remark that Rhoades \cite[Theorem 1.3]{rhoades10} proved the {\itshape cyclic sieving phenomenon} for promotion on standard rectangular tableaux by relating it to cyclic shifting on the Kazhdan--Lusztig basis of the group algebra of the symmetric group $\mathfrak{S}_n$. Purbhoo \cite{purbhoo13} gave another proof of the cyclic sieving phenomenon, by studying the cyclic action of order $k(n-k)$ on $\Gr_{k,n}(\mathbb{C})$ induced by the map $(v_1, \cd, v_n) \mapsto (e^{-i\theta}v_1, e^{-2i\theta}v_2, \cd, e^{-n\theta}v_n)$ on $\mathbb{C}^n$, where $\theta := \frac{2\pi}{k(n-k)}$. It does not appear that his cyclic action is directly related to $\sigma$. Rhoades \cite[Theorem 1.4]{rhoades10} also demonstrated a cyclic sieving phenomenon for promotion on semistandard rectangular tableaux. Shen and Weng \cite{shen_weng} gave an alternative proof of a very similar result, using the cyclic shift map $\sigma$ and its tropicalization. They also explicitly related $\sigma$ to the superpotential $\mathcal{F}_q$ in proving {\itshape cluster duality} for the Grassmannian.

\subsection{Birational rowmotion}\label{sec_rowmotion}

A similar result to Frieden's (discussed in \cref{sec_promotion}) was proved by Grinberg and Roby \cite{grinberg_roby15} in their study of {\itshape birational rowmotion}. This is a map on the torus $(\mathbb{C}^\times)^{k\times (n-k)}$ introduced by Einstein and Propp \cite{einstein_propp14}, which tropicalizes to the {\itshape rowmotion} operation on the poset $[k]\times [n-k]$ (a product of two chains).

There is a connection to the superpotential $\mathcal{F}_q$, as follows. Label each element $(r,s)$ of the poset $[k]\times [n-k]$ by the variable $x_{r,s}$, and add a new minimum element labeled by $1$ and a new maximum element labeled by $q$. Define $\mathcal{L}_q\in\mathbb{Z}[q,x_{r,s}^{\pm 1} : 1 \le r \le k, 1 \le s \le n-k]$ as the sum over all covering relations $u\lessdot v$ of the enlarged poset of the label of $v$ divided by the label of $u$. For example, for $(k,n) = (2,4)$, the labeled poset is
\begin{align}\label{laurent_superpotential}
\begin{tikzpicture}[baseline=(current bounding box.center)]
\tikzstyle{out1}=[inner sep=0,minimum size=1.2mm,circle,draw=black,fill=black,semithick]
\pgfmathsetmacro{\s}{0.84};
\pgfmathsetmacro{\t}{0.06};
\pgfmathsetmacro{\d}{0.72};
\pgfmathsetmacro{\v}{0.72};
\useasboundingbox(-\d-\s,-\t)rectangle(\d+\s,\d+\d+\v+\v+\t);
\node[out1](min)at(0,0){};
\node[out1](11)at($(min)+(0,\v)$){};
\node[out1](12)at($(11)+(\d,\d)$){};
\node[out1](21)at($(11)+(-\d,\d)$){};
\node[out1](22)at($(12)+(-\d,\d)$){};
\node[out1](max)at($(22)+(0,\v)$){};
\node[inner sep=0]at($(min)+(-0.18,0)$){$1$};
\node[inner sep=0]at($(11)+(-0.42,-0.04)$){$x_{1,1}$};
\node[inner sep=0]at($(12)+(0.44,-0.04)$){$x_{1,2}$};
\node[inner sep=0]at($(21)+(-0.40,-0.04)$){$x_{2,1}$};
\node[inner sep=0]at($(22)+(-0.42,0.04)$){$x_{2,2}$};
\node[inner sep=0]at($(max)+(-0.20,0)$){$q$};
\path[semithick](min)edge(11) (11)edge(12) edge(21) (12)edge(22) (21)edge(22) (22)edge(max);
\end{tikzpicture}, \quad\text{ with }\quad \mathcal{L}_q = \frac{x_{1,1}}{1} + \frac{x_{2,1}}{x_{1,1}} + \frac{x_{1,2}}{x_{1,1}} + \frac{x_{2,2}}{x_{2,1}} + \frac{x_{2,2}}{x_{1,2}} + \frac{q}{x_{2,2}}.
\end{align}
The function $\mathcal{L}_q$ was introduced by Eguchi, Hori, and Xiong \cite[Appendix B]{eguchi_hori_xiong97} in constructing a Landau--Ginzburg model for $\Gr_{k,n}(\mathbb{C})$.  Marsh and Rietsch \cite[Proposition 6.10]{marsh_rietsch} showed that $\mathcal{L}_q$ equals the pullback of $\mathcal{F}_q$ to $(\mathbb{C}^\times)^{k\times (n-k)}$, under the embedding $\iota:(\mathbb{C}^\times)^{k\times (n-k)}\hookrightarrow\Gr_{k,n}(\mathbb{C})$ given by
$$
x_{r,s} = \frac{\Delta_{[1,k-r]\cup [k-r+s+1,k+s]}(\iota(x))}{\Delta_{[1,k-r+1]\cup [k-r+s+1,k+s-1]}(\iota(x))}\quad \text{ for } 1 \le r \le k, 1 \le s \le n-k.
$$
(This gives the expression for $\mathcal{F}_q$ which we used in \cref{superpotential_eg}.)

As observed by Galashin and Pylyavskyy \cite[Section 5]{galashin_pylyavskyy19}, the critical points of $\mathcal{L}_q$ in $(\mathbb{C}^\times)^{k\times (n-k)}$ are precisely the fixed points of birational rowmotion. On the other hand, Grinberg and Roby \cite[Proposition 48]{grinberg_roby15} showed that under an embedding very similar to $\iota$, birational rowmotion on $(\mathbb{C}^\times)^{k\times (n-k)}$ becomes the cyclic shift map on $\Gr_{k,n}(\mathbb{C})$. It seems likely that by carefully putting these pieces together, one could more concretely relate $\sigma_t$ and $\mathcal{L}_q$ via birational rowmotion. However, it is important to note that the image of $\iota$ does not in general contain all the fixed points of $\sigma_t$; it is missing those $V\in\Gr_{k,n}(\mathbb{C})$ such that $\Delta_I(V) = 0$ for some $k$-subset $I$ which is a union of two cyclic intervals. Similarly, in general $\mathcal{L}_q$ has fewer critical points than $\mathcal{F}_q$, as seen in \cref{superpotential_eg}, and discussed in more detail by Rietsch \cite[Section 9]{rietsch06}.

\subsection{Control theory and the topology of \texorpdfstring{$\Gr_{k,n}^{\ge 0}$}{the totally nonnegative Grassmannian}}

In \cite{galashin_karp_lam}, Galashin, Karp, and Lam regard $\sigma\in\gl_n(\mathbb{C})$ as a vector field on $\Gr_{k,n}(\mathbb{C})$ (viewed as a quotient of $\GL_n(\mathbb{C})$), whose integral curves are $t\mapsto\exp(t\sigma)(V)$ for $t\in\mathbb{R}$ and $V\in\Gr_{k,n}(\mathbb{C})$. They show that if $V\in\Gr_{k,n}^{\ge 0}$, then $\exp(t\sigma)(V)$ lies in the interior of $\Gr_{k,n}^{\ge 0}$ for $t > 0$, and converges to $V_0$ as $t \to \infty$. (This gives yet another way to see that $V_0$ is the only totally nonnegative fixed point of $\sigma$.) They thereby construct a homeomorphism from $\Gr_{k,n}^{\ge 0}$ onto a small closed ball centered at $V_0$. It would be very interesting to find an analogous flow for each cell of $\Gr_{k,n}^{\ge 0}$, and prove the conjecture of Postnikov that $\Gr_{k,n}^{\ge 0}$ is a regular CW complex \cite[Conjecture 3.6]{postnikov}.\footnote{{\itshape Note added in revision.} Postnikov's conjecture has since been proved using other methods \cite{galashin_karp_lam_regularity}.}

A similar perspective was presented by Ayala, Kliemann, and San Martin \cite{ayala_kliemann_san_martin04}, who used the language of control theory to give an alternative development in type $A$ of Lusztig's theory of total positivity \cite{lusztig}. In this context, $\exp(t\sigma)$ (for $t > 0$) is contained in the interior of the {\itshape compression semigroup} of $\Gr_{k,n}^{\ge 0}$, and so partitions the real Grassmannian $\Gr_{k,n}(\mathbb{R})$ into {\itshape orbits}, each containing a unique critical point of $\sigma$. The action of $\exp(t\sigma)$ as $t\to\infty$ contracts each orbit onto its critical point. There is a unique dense orbit, called the {\itshape stable manifold}, which contains $\Gr_{k,n}^{\ge 0}$, and whose critical point is the {\itshape attractor} $V_0$.

We conclude by again observing that the vector field $\sigma$ and the superpotential of $\Gr_{k,n}(\mathbb{C})$ have the same critical points, several of which are fixed by the twist map.
\begin{prob}\label{relationship_problem}
Explain the relationship between the cyclic shift map, the twist map, and the superpotential.
\end{prob}
It would be especially interesting to make sense of \cref{relationship_problem} for the proper positroid subvarieties of $\Gr_{k,n}(\mathbb{C})$, where the cyclic shift map is not well defined, the twist map is not periodic, and the superpotential is not yet known.

\bibliographystyle{alpha}
\bibliography{Cyclic_shift_bib}

\end{document}